# Equivalence numérique et équivalence cohomologique pour les variétés abéliennes sur les corps finis

By L. Clozel*

Soient $k$ un corps fini de cardinal $q = p^\alpha$, et $A$ une variété abélienne de dimension $g$ sur $k$. Pour $i \leq g$ soit $\mathcal{Z}^i(A)$ le groupe des cycles de codimension $i$ sur $A$, à coefficients dans $\mathbb{Q}$. Considérons pour l'instant la cohomologie étale $\ell$-adique ($\ell \neq p$) $H^\bullet_{\text{ét}}(\overline{A}, \mathbb{Q}_\ell)$ où $\overline{A} = A_k \times \bar{k}$, $\bar{k}$ étant une clôture algébrique de $k$.

On dispose d'une application "classe de cycle" $\text{cl} : \mathcal{Z}^i(A) \to H^{2i}(\overline{A}, \mathbb{Q}_\ell(i))$. Un cycle est *cohomologiquement équivalent* à 0 si $\text{cl}(Z) = 0$; il est *numériquement équivalent* à 0 si le produit d'intersection $Z \cdot T$ est nul pour tout $T \in \mathcal{Z}^{g-i}(A)$.

Puisque $Z \cdot T = \text{cl}(Z)\,\text{cl}(T)$ (au prix de l'isomorphisme canonique $H^{2g}_{\text{ét}}(\overline{A}, \mathbb{Q}_\ell(g)) \cong \mathbb{Q}_\ell$), l'équivalence cohomologique implique l'équivalence numérique. On peut aussi considérer le groupe $\mathcal{Z}^i(\overline{A})$ des cycles sur $\bar{k}$, auquel s'appliquent les mêmes constructions. Le but de cette note est de démontrer:

THÉORÈME 1.   *Soit $\overline{A}$ une variété abélienne sur $\bar{k}$. Il existe alors un ensemble $\mathcal{L}$ de nombres premiers, de densité (analytique ou naïve) non nulle tel que, si $\ell \in \mathcal{L}$, l'équivalence numérique et l'équivalence cohomologique (pour la cohomologie $\ell$-adique) coïncident sur $\overline{A}$.*

COROLLAIRE 1.   *Soit $A$ une variété abélienne sur $k$. Il existe alors un ensemble $\mathcal{L}$ de nombres premiers, de densité non nulle, tel que si $\ell \in \mathcal{L}$, équivalence numérique et cohomologique coïncident sur $A$.*

En effet $\mathcal{Z}(A) = \mathcal{Z}(\overline{A})^G$ où $G = \text{Gal}(\bar{k}/k)$.

Nous appelons *courbe* sur $k$ une courbe projective, lisse et géométriquement irréductible.

COROLLAIRE 2.   *Si $X = \prod_{t=1}^{N} C_t$ est un produit de courbes, il existe un ensemble $\mathcal{L}$ de nombres premiers, de densité (analytique ou naïve) non nulle tel que, si $\ell \in \mathcal{L}$, équivalence numérique et équivalence cohomologique (pour la cohomologie $\ell$-adique) coïncident sur $X$.*

---

*Membre de l'Institut Universitaire de France (1993–1998)

Nous laissons au lecteur le soin de trouver l'extension maximale du Théorème 1 et de son Corollaire.

*Démonstration du Corollaire.* Nous fixons dorénavant un isomorphisme $\mathbb{Z}_\ell \to \mathbb{Z}_\ell(1)$, et nous le ferons dans tout l'article, à l'exception du §4. Ainsi la classe de cycle prend simplement ses valeurs dans $H^{2i}(\overline{A}, \mathbb{Q}_\ell)$. On peut évidemment étendre les scalaires à $\bar{k}$. Notons $\mathcal{A}^\bullet(V)$ le groupe de Chow d'une variété projective $V$, défini par les cycles à coefficients rationnels modulo équivalence linéaire [1]. Si $C_1 \cong \mathbb{P}^1$ (sur $\bar{k}$) on a d'après Grothendieck $\mathcal{A}^\bullet(X) = \mathcal{A}^\bullet(X') \otimes \mathbb{Q}[T]/(T^2)$, $T$ étant la classe fondamentale de $C_1$. (On a posé $X' = \prod_{\ell=2}^N C_t$.) On en déduit aisément que le Corollaire est vrai pour $X$ s'il l'est pour $X'$. On peut donc supposer que le genre $g_t$ de $C_i$ est $\geq 1$ pour tout $t$. Soit $A_t$ la jacobienne de $C_t$, $f_t : C_t \to A_t$ un plongement induisant un isomorphisme $f_t^* : H^1(A_t) \to H^1(C_t)$ où l'on note simplement $H^i(\cdot) = H^i_{\text{ét}}(\cdot, \mathbb{Q}_\ell)$ et $f : X \to A = \prod A_t$ le produit des $f_t$.

Soit $g = \sum_{t=1}^N g_t$. On a alors deux applications

$$f_* : H^i(X) \to H^{2(g-N)+i}(A),$$
$$f^* : H^i(A) \to H^i(X)$$

duales pour la dualité de Poincaré. Les applications $f_*$ et $f^*$ sont aussi définies entre les groupes de Chow.

L'application $f^*$, en cohomologie, est surjective puisque $H^\bullet(X)$ est engendré comme algèbre par $H^1(X) = \bigoplus H^1(C_t)$ et $f^* : H^1(A) \to H^1(X)$ est bijectif. Donc $f_*$ est injectif. Soit alors $Z \in \mathcal{Z}^i(X)$ et supposons $Z$ numériquement égal à 0. Alors $f_*Z$ l'est aussi d'après la formule d'adjonction. D'après le Théorème 1, $f_*Z$ est alors cohomologiquement trivial, et il en est de même de $Z$. Ceci termine la démonstration.

En caractéristique nulle, le fait que l'équivalence cohomologique et numérique coïncident résulte de l'algébricité de l'opérateur $\Lambda$ de la théorie de Lefschetz et du théorème de positivité de Hodge; le théorème est alors dû à D. Lieberman (cf. [6]). Notre preuve du Théorème 1 repose aussi sur l'algébricité de $\Lambda$; l'usage de la positivité de Hodge est remplacé par l'abondance d'endomorphismes de $A$.

Je remercie S. Bloch, J.-L. Colliot-Thélène, N. Katz, L. Illusie, G. Laumon, M. Raynaud et D. Ramakrishnan pour d'utiles conversations. Je suis tout particulièrement reconnaissant à Ofer Gabber pour avoir décelé une erreur dans une première version de ce travail. Le rapporteur (P. Deligne) m'a indiqué comment interpréter la démonstration à l'aide du groupe d'automorphismes de $H^\bullet(\overline{A}, \mathbb{Q}_\ell)$ engendré par $\text{End}(A)^*$ et les opérateurs de Lefschetz. Ceci sera publié ultérieurement; je le remercie très vivement.

## 1. Généralités sur la multiplication complexe

Dans ce qui suit $A$ est une variété abélienne définie sur $k$. Supposons d'abord $A \times \overline{k}$ irréductible. Soit $\operatorname{End}^0 A = \operatorname{End}_k(A) \otimes \mathbb{Q}$. On sait que $A$ admet toujours des multiplications complexes: il existe un corps $E$ de degré $2g$ sur $\mathbb{Q}$ qui est un corps de multiplication complexe (corps CM), i.e. une extension quadratique totalement imaginaire d'un corps de nombres totalement réel, contenu dans $\operatorname{End}^0 A$. On note $c$ la conjugaison complexe de $E$.

Soit $M$ un corps CM, extension galoisienne de $\mathbb{Q}$, dans lequel se plonge $E$; l'existence d'un tel corps est bien connue [9]. Si $\Pi$ désigne l'ensemble des plongements de $E$ dans $M$, on a alors

$$(1.1) \qquad E \underset{\mathbb{Q}}{\otimes} M \cong \bigoplus_{\sigma \in \Pi} M,$$

l'isomorphisme étant donné par

$$e \otimes m \mapsto (\sigma(e)m)_\sigma.$$

Donc (1.1) est $M$-linéaire pour l'action à droite de $M$ sur $E \otimes M$, et chacune de ses composantes est $\sigma$-linéaire pour l'action à gauche de $E$.

Soit $\ell$ un nombre premier $\neq p$, et soit $\lambda$ une place de $M$ divisant $\ell$. On a de même

$$(1.2) \qquad E \underset{\mathbb{Q}}{\otimes} M_\lambda \cong \bigoplus_\sigma M_\lambda,$$

l'isomorphisme étant donné par $e \otimes m \mapsto (\sigma(e)m)_\sigma$.

On notera simplement $H^i(A, M_\lambda)$ le module $H^i_{\text{ét}}(\overline{A}, \mathbb{Q}_\ell) \otimes_{\mathbb{Q}_\ell} M_\lambda$. Rappelons que $H^1_{\text{ét}}(\overline{A}, \mathbb{Q}_\ell)$ est un module libre de rang 1 sur $E \otimes \mathbb{Q}_\ell$, $E$ opérant par multiplications complexes. On en déduit que $H^1(A, M_\lambda)$ est libre de rang 1 sur $E \underset{\mathbb{Q}}{\otimes} M_\lambda = (E \otimes \mathbb{Q}_\ell) \underset{\mathbb{Q}_\ell}{\otimes} M_\lambda = \bigoplus_\sigma M_\lambda$, donc libre de rang $2g$ sur $M_\lambda$. On a alors:

LEMME 1.1.  $H^1(A, M_\lambda) = \bigoplus_{\sigma \in \Pi} H^1(A, M_\lambda)_\sigma$, chaque facteur $H^1(A, M_\lambda)_\sigma$ étant l'image de $H^1(A, M_\lambda)$ par un projecteur algébrique.

En effet, soit $e_\sigma \in \bigoplus_\tau M$ l'élément de composante 1 à la place $\sigma$ et 0 en $\tau \neq \sigma$: alors $e_\sigma$ définit sous (1.1) un élément de $\operatorname{End}^0(A) \otimes M$ qui projette sur $H^1(A, M_\lambda)_\sigma$. On note ce projecteur $p_\sigma$.

Si $\lambda \in E$, on désigne par $[\lambda]$ l'endomorphisme associé de $A$ et son effet sur les espaces de cohomologie.

LEMME 1.2.  $H^1(A, M_\lambda)_\sigma = \{\omega \in H^1(A, M_\lambda) : [\lambda]\omega \equiv \lambda^\sigma \omega\}$; il est de rang 1 sur $M_\lambda$.

C'est clair.

On va maintenant décomposer la cohomologie en degré $i$ sous l'action de $E$. Rappelons que $c$ désigne la conjugaison complexe de $E$. Soit $\Sigma$ un type CM de $E$, i.e., un ensemble de plongements de $E$ dans $M$ tel que tout plongement $E \to M$ soit de la forme $\sigma$ ou $c\sigma$ pour un $\sigma \in \Sigma$. On sait que $H^i_{\text{ét}}(\overline{A}, \mathbb{Q}_\ell) = \mathbf{\Lambda}^i H^1_{\text{ét}}(\overline{A}, \mathbb{Q}_\ell)$, le produit extérieur étant pris sur $\mathbb{Q}_\ell$. On a donc de même $H^i(A, M_\lambda) = \mathbf{\Lambda}^i H^1(A, M_\lambda)$ (produit extérieur sur $M_\lambda$). On notera $I, J$ des sous-ensembles de $\Sigma$ et $^c\Sigma$ respectivement; si $\lambda \in E$, soit $\lambda^I \lambda^J = \prod_{\sigma \in I \amalg J} \lambda^\sigma \in M$.

LEMME 1.3. $\quad H^i(A, M_\lambda) = \bigoplus_{\substack{I,J \\ |I|+|J|=i}} H^i(A, M_\lambda)_{IJ}$ où
$$H^i(A, M_\lambda)_{IJ} = \{\omega \in H^i : [\lambda]\omega \equiv \lambda^I \lambda^J \omega\}$$
est de rang 1 sur $M_\lambda$.

Ceci résulte immédiatement du Lemme 1.2. Notons aussi qu'on a des projecteurs $p_{IJ}$ sur $H^i(A, M_\lambda)_{IJ}$ définis comme après le Lemme 1.

Puisque $E$ est une sous-algèbre abélienne maximale de $\text{End}^0(A)$, l'endomorphisme de Frobenius $\pi$ peut être considéré comme un élément de $E$. Pour tout $\sigma \in \Sigma$, $\pi^\sigma \pi^{c\sigma} = \pi^\sigma \bar{\pi}^\sigma = q$. Par conséquence $\pi$ opère sur $H^2(A, M_\lambda)_{\sigma,c\sigma}$ comme multiplication par $q$.

Notons que l'application cycle envoie $\mathcal{Z}^i(A) \otimes M$ vers $H^{2i}(\overline{A}, \mathbb{Q}_\ell) \otimes M_\lambda = H^{2i}(A, M_\lambda)$.

LEMME 1.4. $\quad$ Pour tout $\sigma \in \Sigma$, il existe un cycle $L_\sigma \in \mathcal{Z}^1(A) \otimes M$ tel que $\text{cl}(L_\sigma)$ soit un générateur de $H^2(A, M_\lambda)_{\sigma,c\sigma}$ sur $M_\lambda$.

*Démonstration.* L'élément de Frobenius $\pi$ opère sur $H^2(A, M)_{\sigma,c\sigma}$ par la valeur propre $q$. Le Lemme se déduit alors aisément, par extension des scalaires, du résultat analogue de Tate pour $H^2(\overline{A}, \mathbb{Q}_\ell)$ [11].

*Remarque.* Dans le §2, le Lemme 1.4 sera retrouvé (au moins pour un choix particulier de $\Sigma$) à l'aide du relèvement à la caractéristique nulle. Il n'est donc pas utilisé *directement* dans nos démonstrations (mais la théorie de Honda-Tate repose sur les résultats de Tate [11], donc cette indépendance est vraisemblablement illusoire).

Nous considérons maintenant l'action de $\text{Gal}(M/\mathbb{Q})$ sur la cohomologie. Si $\tau \in \text{Gal}(M/\mathbb{Q})$, $\tau$ définit par transport de structure un isomorphisme de $H^i(A, M_\lambda)$ avec $H^i(A, M_{\tau\lambda})$ où $\tau\lambda$ est la place de $M$ déduite de $\lambda$ par l'action de $\tau$.

LEMME 1.5. *Soit $c \in \mathrm{Gal}(M/\mathbb{Q})$ la conjugaison complexe. Alors $c$ envoie $H^i(A, M_\lambda)_{I,J}$ vers $H^i(A, M_{c\lambda})_{cJ,cI}$.*

*Démonstration.* Soit $\omega \in H^i(A, M_\lambda)$ telle que $[\lambda]\omega = \lambda^I \lambda^J \omega$. Alors $^c([\lambda]\omega) = [\lambda]^c\omega$ puisque $[\lambda]$ opère par un endomorphisme de $\overline{A}$. On a donc $[\lambda]\,^c\omega = {}^c(\lambda^I \lambda^J \omega) = {}^c\lambda^{Ic} \lambda^{J\,c}\omega$ (d'après (1.2)) $= \lambda^{cJ} \lambda^{cI\,c}\omega$. L'assertion résulte du Lemme 1.2.

COROLLAIRE 1.6. *Si $H^{2i}(A, M)_{I,J}$ est engendré sur $M_\lambda$ par la classe d'un cycle algébrique dans $\mathcal{Z}^i(A) \otimes M$, il en est de même de $H^{2i}(A, M_{c\lambda})_{cJ,cI}$.*

En effet $c$ opère (sur les coefficients) en préservant les cycles algébriques: on est ramené au Lemme 1.5.

## 2. Division des cycles algébriques

On va maintenant énoncer un résultat simple, mais fondamental, de division des cycles algébriques. Quitte à remplacer $k$ par une extension finie, on sait [12] que $A$, ou au moins une variété $A'$ sur $k$ isogène à $A$, admet un relèvement en caractéristique 0: il existe un corps de nombres $F$ et une variété $\widetilde{A}$ sur $F$ ayant les propriétés suivantes: il existe une place $v$ de $F$, de corps résiduel $k$, telle que le modèle de Néron $\mathcal{A}$ de $\widetilde{A}$ sur $\mathcal{O}_F$ ait bonne réduction en $v$ et que $\mathcal{A} \times_{\mathcal{O}_F} k \cong A'$. L'action de $E$ sur (la variété abélienne à isogénie près déduite de) $A'$ se relève à $\widetilde{A}$, et cette action est compatible avec la réduction en $v$.

Noter que le Théorème 1 est invariant par isogénie; puisque l'action de $G = \mathrm{Gal}(\overline{k}/k)$ sur $\mathcal{Z}(\overline{A})$ est localement finie, il suffit pour démontrer le théorème de vérifier son corollaire pour toute extension assez grande $k'$ du corps $k$. Nous supposons donc désormais que $A$ est définie sur $k$ et admet un relèvement en caractéristique 0, au sens précisé ci-devant. L'action de $E$ sur l'espace tangent de $\widetilde{A}$ en 0 définit alors un type CM, que l'on supposera égal à $\Sigma$.

Notons $F_v$ la complétion de $F$ en $v$, $\mathcal{O}_v$ son anneau d'entiers et $\mathcal{A} = \mathcal{A} \times_{\mathcal{O}_F} \mathcal{O}_v$: $\mathcal{A}_v$ est un schéma abélien sur $\mathcal{O}_v$. On sait que $\mathcal{A}_v$ admet un plongement projectif $j : \mathcal{A} \to \mathbb{P}^N_{\mathcal{O}_v}$ correspondant à un diviseur très ample $\mathcal{L} = j^*\mathcal{O}(1)$ sur $\mathcal{A}_v$.

Quitte à remplacer $k$ par une extension finie, on peut trouver, d'après le théorème de Bertini, une section de $j^*\mathcal{L}$ sur $\mathcal{A}_v \times k$ dont le lieu des zéros est une sous-variété lisse de $\mathcal{A}_v \times k$. D'après un théorème de Mumford [8], $H^1(\mathcal{A}_v, \mathcal{L}) = \{0\}$, et donc cette section se relève en un élément de $H^0(\mathcal{A}_v, \mathcal{L})$. En considérant le lieu de ses zéros, on obtient donc un sous-schéma lisse $L_v$ de $\mathcal{A}_v$ tel que $L_v \times k$ est une section hyperplane de $\mathcal{A}_v \times k$ et $L_v \times F'_v$ une section hyperplane de $\mathcal{A}_v \times F'_v$ ($F'_v$ étant une extension de $F_v$ de corps résiduel $k$). Quitte à remplacer aussi $F$ par une extension finie, ce que nous oublions désormais dans la notation, on voit qu'on a une classe de cohomologie $L$ dans $H^2(\widetilde{A} \times \overline{F}, \mathbb{Q}_\ell) \cong H^2(\overline{A}, \mathbb{Q}_\ell)$ qui est, à la fois sur $\widetilde{A} \times \overline{F}$ et sur $\overline{A}$, la classe de cohomologie d'une section hyperplane.

LEMME 2.1. *Les composantes isotypiques de $L$ dans la décomposition $H^2(A, M_\lambda) = \bigoplus_{I,J} H^2(A, M_\lambda)_{I,J}$ sont de type $(\sigma, c\sigma)$, $\sigma$ parcourant $\Sigma$. Chaque composante de type $(\sigma, c\sigma)$ est non nulle.*

*Démonstration.* Il suffit de le vérifier sur $\widetilde{A}$, et même sur $B = \widetilde{A} \times_F \mathbb{C}$ pour un plongement complexe de $F$. Choisissons un plongement complexe de $M$. A isogénie près (ce qui ne change rien au problème) $B$ s'identifie à $\mathbb{C}^g/\Sigma(\mathcal{O})$, où $\mathcal{O}$ est l'anneau des entiers de $E$ et où $\Sigma$ envoie $\mathcal{O}$ dans $\mathbb{C}^g$ par composition avec l'injection de $M$ dans $\mathbb{C}$. En particulier $H_1(B) \otimes \mathbb{Q} = \mathcal{O} \otimes \mathbb{Q} = E$. Noter que $\widetilde{A}$, et donc $B$, est absolument irréductible. On sait alors, sous nos hypothèses, que $\mathrm{End}^0(B) \cong E$ (Shimura-Taniyama; cf. Shimura [9, p. 128]). D'après Shimura (ibid.), la classe de cohomologie associée à une polarisation est alors donnée (en cohomologie de Betti) par une forme de Riemann sur $E$:

$$(x, y) \mapsto \mathrm{Tr}_{E/\mathbb{Q}}(\zeta x^c y)$$

où $\zeta \in E$ vérifie $^c\zeta = -\zeta$.

Par extension des scalaires, on obtient une forme alternée sur $E \otimes M \cong \bigoplus_\sigma M$ donnée par

$$((x_\sigma), (y_\sigma)) \mapsto \sum_\sigma x_\sigma {}^\sigma\zeta \, y_{c\sigma}.$$

Son type (et donc celui de $L$) est bien celui indiqué par le lemme; il est clair que chaque composante est non nulle.

Soit $(I, J)$ un type, avec $I \subset \Sigma$ et $J \subset {}^c\Sigma$. Soit $K = I \cap {}^cJ$, $I^0 = I - K$, $J^0 = J - {}^cK$. On a donc $I = I^0 \amalg K$, $J = J^0 \amalg {}^cK$. La démonstration du Théorème 1 va reposer sur le résultat suivant.

PROPOSITION 2.2. *Soit $(I, J)$ un type, $|I| + |J| = 2j$, et supposons qu'il existe un cycle $Z \in \mathcal{Z}^j(A) \otimes M$ dont la classe engendre $H^{2j}(A, M_\lambda)_{I,J}$. Alors il existe un cycle $T \in \mathcal{Z}^{j-k}(A) \otimes M$ (avec $k = |K| = |I \cap {}^cJ|$) dont la classe engendre $H^{2j-2k}(A, M_\lambda)_{I^0, J^0}$.*

*Démonstration.* On notera simplement $L : H^i \to H^{i+2}$ (où $H^i = H^i(A, M_\lambda)$) la multiplication par la classe de Lefschetz. On considère l'opérateur $\Lambda : H^i \to H^{i-2}$ de la théorie de Lefschetz. On sait qu'il est algébrique, d'après un théorème de Liebermann [6]; d'après Katz et Messing [5], les projections de $H^\bullet$ sur ses composantes $H^i$ sont algébriques; enfin, pour $i$ fixé, la décomposition $H^i = \bigoplus_{I,J} H^i_{I,J}$ est algébrique. Il en résulte que pour tout $I, J$ et $I', J'$ (avec $|I| + |J| = |I'| + |J'| + 2$) la composante de $\Lambda$ qui envoie $H^\bullet_{I,J}$ vers $H^\bullet_{I',J'}$ est algébrique.

Fixons une base $(\omega_\sigma)$ ($\sigma \in \Sigma \cup {}^c\Sigma$) de $H^1(A, M_\lambda)$ sur $M_\lambda$; d'après le Lemme 1.7, on peut supposer que $\omega_\sigma \omega_{c\sigma} = L_\sigma$ est une composante isotypique de la classe de Lefschetz, donc algébrique. On a $L = \sum_{\sigma \in \Sigma} L_\sigma$, et $L_\sigma^2 = 0$. Pour tout $i \leq g$, $L^{g-i} : H^i \to H^{2g-i}$ a donc pour expression

$$(2.1) \qquad L^{g-i} = \left(\sum_\sigma L_\sigma\right)^{g-i} = (g-i)! \sum_K L_K$$

où la somme porte sur les sous-ensembles $K$ de $\Sigma$ de cardinal $(g-i)$, et

$$(2.2) \qquad L_K = \prod_{\sigma \in K} L_\sigma .$$

Si $(I, J)$ est un couple d'indices, soit $\omega_{IJ} = \prod_{\sigma \in I \amalg J} \omega_\sigma$. D'après (2.1),

$$(2.3) \qquad L^{g-i} \omega_{IJ} = \sum_K (g-i)!\, \omega_{I \amalg K, J \amalg {}^c K}$$

où la somme porte sur les $K$ (avec $|K| = g - i$) disjoints de $I$ et $J$.

Soit $\theta_i : H^{2g-i} \to H^i$ l'inverse de $L^{g-i}$. On a alors la conséquence suivante:

LEMME 2.3.  *Supposons que $|I'| + |J'| = 2g - i$. Si $\theta_i(\omega_{I',J'})$ a une composante non nulle en $(I, J)$ – avec $|I| + |J| = i$ – alors $I' = I \amalg K$, $J' = J \amalg {}^c K$ avec $K \subset \Sigma$ et $|K| = g - i$.*

Supposons $2 \leq i \leq g$. Alors $\Lambda : H^i \to H^{i-2}$ est défini (cf. Kleiman [7, p. 13]) par

$$\Lambda = \theta^{i-2} L^{g-i+1} .$$

On déduit alors de (2.3) et du Lemme 2.3:

LEMME 2.4.  *Soit $(I, J)$ avec $|I| + |J| = i$, et supposons que $\Lambda \omega_{IJ}$ a une composante non nulle en $(I', J')$ avec $|I'| + |J'| = i - 2$. Alors il existe $K, L \subset \Sigma$ avec $|K| = g - i + 1$, $|L| = g - i + 2$ tels que*

$$(2.4) \qquad I \amalg K = I' \amalg L ,$$

$$(2.5) \qquad J \amalg {}^c K = J' \amalg {}^c L .$$

(En particulier, $K$ est disjoint de $I \cup {}^c J$ et $L$ de $I' \cup {}^c J'$).

Soit $I_1 = I \amalg K$, $J_1 = J \amalg {}^c K$ et définissons $I^0$, $I_1^0$, $J^0$, $J_1^0$ comme avant la Proposition 1.8. On vérifie aussitôt que $I^0 = I_1^0$, $J^0 = J_1^0$. Il résulte alors de (2.4) et (2.5) que $I^0 = (I')^0$, $J^0 = (J')^0$. Supposons alors qu'il existe un cycle algébrique de type $(I, J)$. Si $\Lambda\omega_{IJ} \neq 0$, le Lemme 2.4 implique alors l'existence d'un cycle algébrique de type $(I', J')$ avec $|I'|+|J'| = |I|+|J|-2$ et $I^0 = (I')^0$, $J^0 = (J')^0$. Pour $i = 2j \leq g$, la Proposition 2.2 résulte alors, par récurrence sur le degré, du lemme suivant:

LEMME 2.5. *Si $K = I \cap {}^c J \neq \emptyset$, une classe $\omega$ de type $(I, J)$ n'est pas primitive, et donc $\Lambda\omega \neq 0$.*

*Démonstration.* D'après (2.3), on a

$$(2.6) \qquad L^{g-i+1}\omega_{IJ} = \sum_L (g-i+1)!\, \omega_{I \amalg L, J \amalg {}^c L}$$

où $|L| = g - i + 1$. Si $I \cap {}^c J \neq \emptyset$, $|I \cup {}^c J| \leq i - 1$ où $i = |I| + |J|$. Il existe donc $g - i + 1$ éléments de $\Sigma$ n'appartenant pas à $I \cup {}^c J$: si $L$ est formé de ces éléments, le terme $\omega_{I \amalg L, J \amalg {}^c L}$ de (2.6) est non nul, donc $\omega_{IJ}$ n'est pas primitive.

Il reste à traiter le cas où $i = 2j > g$. L'opérateur $\Lambda$ est alors défini par le diagramme commutatif [7, p. 13]:

$$(2.7) \qquad \begin{array}{ccc} H^{2g-i} & \xrightarrow{L^{i-g}} & H^i \\ L \downarrow & & \downarrow \Lambda \\ H^{2g-i+2} & \xrightarrow{L^{i-2-g}} & H^{i-2} \end{array},$$

i.e., $\Lambda = L^{i-1-g}\theta_{2g-i}$. (Si $i = g + 1$, $\Lambda$ est simplement $\theta_{g-1}$; $\theta_{2g-i}$ est l'inverse de $L^{i-g} : H^{2g-i} \to H^i$). On en déduit l'analogue du Lemme 2.4: si $|I|+|J| = i$ et si $\Lambda\omega_{IJ}$ a une composante non nulle en $(I', J')$ avec $|I'|+|J'| = i-2$, alors

$$(2.8) \qquad I = I'' \amalg K\ ,\ J = J'' \amalg {}^c K\ ,$$

$$(2.9) \qquad I = I'' \amalg L\ ,\ J' = J'' \amalg {}^c L$$

avec $|K| = i - g$, $|L| = i - g - 1$. On a de nouveau $I^0 = (I'')^0$, $J^0 = (J'')^0$ d'où $I^0 = (I')^0$, $J^0 = (J')^0$. De plus, une classe de type $(I, J)$ où $|I|+|J| > g$ ne peut être primitive (c'est du reste clair sur (2.7) puisque $L^{i-1-g} : H^{2g-i} \to H^{i-2}$ est injectif), et l'on termine comme dans le cas précédent.

On notera la conséquence suivante de la Proposition 2.2, que nous n'avons obtenue que de façon indirecte:

COROLLAIRE 2.6. *Supposons qu'il existe un cycle $Z \in \mathcal{Z}^j(A) \otimes M_\lambda$ dont la classe est non nulle et de type $(I, J)$. Alors, pour tout $K \subset I \cap {}^c J$, il existe un cycle algébrique dont la classe engendre la cohomologie de type $(I-K, J-{}^c K)$.*

En d'autres termes, on peut *diviser* les cycles algébriques isotypiques par les composantes isotypiques de la classe de Lefschetz. Le Corollaire se déduit de la Proposition 2.2 en descendant à $(I^0, J^0)$ puis en *multipliant* par des composantes isotypiques de la classe de Lefschetz.

## 3. Le cas réductible; démonstration du Théorème 1

Une variété abélienne arbitraire $A/\bar{k}$ est isogène à un produit de composantes irréductibles ($\equiv$ simples) $A_t$ ($t = 1, \ldots N$). On peut supposer que chaque $A_t$ se relève en caractéristique nulle. Soit $E_t$ ($t = 1, \ldots N$) un corps de multiplication complexe pour $A_t$. Alors $E = \bigoplus_{t=1}^{N} E_t$ opère sur $A$. Soit $M$ le composé des clôtures galoisiennes $M_t$ des $E_t$. Alors on a une décomposition isotypique

$$(3.1) \qquad H^\bullet(A, M_\lambda) = \bigoplus_{\substack{I=(I_1,\ldots I_N) \\ J=(J_1,\ldots J_N)}} H^\bullet(A, M_\lambda)_{IJ}$$

où l'on a choisi des types CM $\Sigma_t$ pour les corps $E_t$, et où $I_t \subset \Sigma_t$ et $J_t \subset {}^c\Sigma_t$.

Les composantes de la classe de Lefschetz fournissent des cycles algébriques de type $(\sigma, c\sigma)$ où $\sigma$ parcourt $\Sigma_1 \amalg \cdots \amalg \Sigma_N$.

Définissons alors $(I^0, J^0)$ par ces composantes $(I_t^0, J_t^0)$, elles-mêmes définies comme avant la Proposition 2.2. Soit $k = \sum_t |K_t|$, la notation étant évidente. On peut alors "künnethiser" les démonstrations du §3, en remplaçant partout l'action d'un corps par l'action de l'algèbre $E = \bigoplus E_t$. On obtient ainsi l'analogue de la Proposition 2.2:

LEMME 3.1.    *Soit $(I, J)$ un type, avec $|I| + |J| = 2j$, et supposons qu'il existe un cycle $Z \in \mathcal{Z}^j(A) \otimes M$ dont la classe engendre $H^{2j}(A, M_\lambda)_{I,J}$. Alors il existe $T \in \mathcal{Z}^{j-k}(A) \otimes M$ dont la classe engendre $H^{2j-2k}(A, M_\lambda)_{I^0, J^0}$.*

Nous pouvons maintenant démontrer le Théorème 1. Soit $A/k$ une variété abélienne de dimension $g$; nous supposons que $A$ est un produit de variétés simples se relevant en caractéristique 0, ce qui est loisible à isogénie près.

Rappelons la construction des opérateurs $p_{IJ}$ dans $H^{2j}(A, M_\lambda)$. Pour $I, J = (I_t), (J_t)$ comme ci-dessus, soit $e_{IJ}$ l'élément de $E \otimes M = \bigoplus E_t \otimes M$ de composantes 1 aux places $\sigma \in I_t \cup J_t$, 0 ailleurs. On dispose alors de l'endomorphisme $p_{IJ} = [e_{IJ}]$ de $H^i$; c'est un projecteur et *en degré $i = |I|+|J|$*, on vérifie aisément que $p_{IJ}$ est le projecteur sur la composante de type $(I, J)$ de (3.1).

Si $Z \in \mathcal{Z}^j(A) \otimes M$, $[e_{IJ}]Z$ est un cycle algébrique et l'on a

$$(3.2) \qquad \mathrm{cl}([e_{IJ}]Z) = p_{IJ}\,\mathrm{cl}(Z)$$

et donc

$$\text{(3.3)} \qquad \text{cl}(Z) = \sum_{I,J} \text{cl}([e_{IJ}]Z) \, .$$

LEMME 3.2.   *Si $Z \in \mathcal{Z}^j(A) \otimes M$ est numériquement équivalent à 0, $[e_{IJ}]Z$ l'est aussi pour tout $I, J$.*

En effet pour tout $T \in \mathcal{Z}^{g-j}(A) \otimes M$, on a alors $[e_{IJ}]Z \cdot T = Z \cdot {}^t[e_{IJ}]T$, la transposée d'une correspondance étant définie comme dans [6]. Mais $[e_{IJ}]$ et donc sa transposée sont des correspondances algébriques.

Supposons alors $Z$, et donc $[e_{IJ}]Z$, numériquement égal à 0. Soit $\omega = \text{cl}(Z)$, $\omega_{IJ} = p_{IJ}\omega = \text{cl}([e_{IJ}]Z)$, et supposons $\omega_{IJ} \neq 0$; c'est alors la classe d'un cycle. D'après le Lemme 3.1, il existe un cycle $T$ dans $\mathcal{Z}^{j-k}(A) \otimes M$, tel que $\text{cl}(T)$ est de type $(I^0, J^0)$ et que

$$\text{(3.4)} \qquad \omega_{IJ} = \mu \, L_K \, \text{cl}(T) \, , \ \mu \in M_\lambda \, , \ \mu \neq 0 \, .$$

(On a désigné par $L_K = \prod_t L_{K_t}$ le produit convenable des composantes de la classe de Lefschetz; notation (2.2).)

Considérons alors la conjugaison complexe $c \in \text{Gal}(M/\mathbb{Q})$. D'après l'extension évidente du Lemme 1.5 et de son Corollaire, ${}^c(\text{cl}(T))$ est de type ${}^cJ^0$, ${}^cI^0$ et évidemment non nul; c'est une classe de cohomologie dans $H^\bullet(A, M_{c\lambda})$. Supposons alors $\lambda$ fixe par la conjugaison complexe. Puisque $I^0 \cap {}^cJ^0 = J^0 \cap {}^cI^0 = \emptyset$, on a nécessairement $\text{cl}(T) \cdot {}^c\text{cl}(T) \neq 0$ puisque l'algèbre de cohomologie est libre sur les $\omega_{IJ}$. Puisque $K = \coprod K_t$ est disjoint de $I^0 \cup {}^cJ^0$, on en déduit de même que $L_K \text{cl}(T)^c \text{cl}(T) \neq 0$; d'après (3.4), $\omega_{IJ} \text{cl}({}^cT) \neq 0$. C'est une classe de cohomologie algébrique de type $(I \amalg {}^cJ^0, J \amalg {}^cI^0 = (K \amalg I^0 \amalg {}^cJ^0, {}^cK \amalg J^0 \amalg {}^cI^0)$. Si $H$ est le complément dans $\Sigma = \coprod \Sigma_t$ de $K \amalg I^0 \amalg {}^cJ^0$, on en déduit enfin que $L_H \omega_{IJ} \text{cl}({}^cT) \neq 0$. Mais ceci contredit le Lemme 3.2 puisque $L_H$, et ${}^cT$, sont algébriques. Donc $\omega_{IJ} = 0$ pour tous $I, J$ et ceci démontre le Théorème 1.

Il reste à vérifier que l'on peut bien trouver une famille de $\ell$, de densité non nulle, tels qu'il existe une place $\lambda | \ell$ de $M$ telle que $c\lambda = \lambda$. Soit $R$ le sous-corps réel maximal de $M$. (Noter que $M$ est un corps CM: [9, Prop. 5.12].) On a une suite exacte

$$\text{(3.5)} \qquad 1 \to \{1, c\} \to \text{Gal}(M/\mathbb{Q}) \to \text{Gal}(R/\mathbb{Q}) \to 1$$

et l'élément $c$ est central dans $\text{Gal}(M/\mathbb{Q})$. Il suffit pour notre démonstration de supposer que $c$ appartient au groupe de décomposition $D_\lambda$ de $\lambda$. Si $M$ est non ramifié en $\ell$, ceci revient à supposer que $c$ est une puissance d'un élément de Frobenius $\mathfrak{F} \in \text{Gal}(M/\mathbb{Q})$, ce qui est vrai par exemple si $\mathfrak{F} = c$; d'après le théorème de Čebotarev, ceci est vrai avec la fréquence indiquée par le Théorème 1.

## 4. Relation avec d'autres conjectures

D'après le Théorème 1, on pourrait démontrer (sans condition sur $\ell$) que l'équivalence cohomologique et numérique coïncident sur $A$ si l'on disposait d'une quelconque des propriétés suivantes. La première est bien connue:

CONJECTURE 4.1. *Si $Z \in \mathcal{Z}^i(A)$ et si $\text{cl}_\ell(Z) \in H^{2i}(\overline{A}, \mathbb{Q}_\ell(i)) \neq 0$ pour un $\ell \neq p$, alors $\text{cl}_m(Z) \in H^{2i}(\overline{A}, \mathbb{Q}_m(i)) \neq 0$ pour tout nombre premier $m \neq p$.*

En effet, on pourrait alors choisir $m$ vérifiant les hypothèses du Théorème, et en déduire que $Z$ n'est pas numériquement équivalent à 0.

Notons qu'une variante *à caractéristique* (résiduelle) *constante* de la Conjecture 4.1 ferait aussi l'affaire. Soient $M$ un corps de nombres et $Z$ un cycle *à coefficients dans $M$*. Soient $\lambda, \mu$ deux places de $M$ divisant $\ell \neq p$.

CONJECTURE 4.2. *Si $Z \in \mathcal{Z}^i(A, M)$ alors $\text{cl}_\lambda(Z) \in H^{2i}(A, M_\lambda(i)) = 0$ si et seulement si $\text{cl}_\mu(Z) \in H^{2i}(A, M_\mu(i)) = 0$.*

Je ne vois pas, cependant, pourquoi cette conjecture serait plus accessible que la Conjecture 4.1.

Comme me l'a fait remarquer D. Ramakrishnan, la démonstration donnée ici s'adapterait si l'on disposait d'une propriété conjecturée par C. Soulé.

Soit $CH^i(A, \mathbb{Q})$ le groupe des classes de cycles modulo équivalence linéaire, à coefficients dans $\mathbb{Q}$; soit $CH^i(A, \mathbb{Q}_\ell) = CH^i(A, \mathbb{Q}) \otimes \mathbb{Q}_\ell$.

CONJECTURE 4.3 (Soulé). *L'application classe de cycle donne un isomorphisme*

$$(4.1) \qquad CH^i(A, \mathbb{Q}_\ell) \to H^{2i}(\overline{A}, \mathbb{Q}_\ell(i))^G$$

*où $G = \text{Gal}(\bar{k}/k)$.*

Pour cette conjecture voir Soulé [10]: elle n'est pas formulée explicitement dans cet article mais une forme plus faible y est *démontrée* pour $i = 0, 1, d-1$ où $d$ est la dimension de $A$.

PROPOSITION 4.4. *Si l'application (4.1) est injective, l'équivalence numérique et cohomologique coïncident pour $A$ et $\ell$.*

Le point notable est que nous ne supposons pas la surjectivité, et donc la conjecture de Tate ([13: Conj. $T^j$]).

*Démonstration.* Soit $M_\ell = M \otimes \mathbb{Q}_\ell$, avec les notations habituelles. Nous négligeons aussi les twists à la Tate dans les notations. On a par hypothèse

$$(4.2) \qquad CH^i(A, \mathbb{Q}) \otimes \mathbb{Q}_\ell \hookrightarrow H^{2i}(\overline{A}, \mathbb{Q}_\ell)$$

d'où, en prenant le produit tensoriel (plat) par $M$:

$$(4.3) \qquad CH^i(A, M) \otimes \mathbb{Q}_\ell \hookrightarrow H^{2i}(\overline{A}, M_\ell)$$

et en appliquant les projecteurs $p_{IJ}$:

$$(4.4) \qquad (CH^i(A, M)_{IJ}) \otimes \mathbb{Q}_\ell \hookrightarrow H^{2i}(\overline{A}, M_\ell)_{IJ}$$

– noter que $CH^i(A, M)$ se décompose selon les projecteurs $p_{IJ}$ vu l'injectivité de (4.3). Supposons alors qu'il existe un cycle de type $(I, J)$ dont la classe de cohomologie dans $H^{2i}(\overline{A}, M_\ell)$ est non nulle. Puisque le membre de droite de (4.4) est de rang 1 sur $M_\ell$, $CH^i(A, M)_{IJ}$ est alors de rang 1 sur $M$, et l'application (4.4) est nécessairement un isomorphisme sur $M_\ell = M \otimes \mathbb{Q}_\ell$. Si $Z_{IJ}$ est un tel cycle, il en résulte que toutes ses composantes dans $H^{2i}(\overline{A}, M_\ell)_{IJ} = \bigoplus_{\lambda | \ell} H^{2i}(A, M_\ell)_{IJ}$ sont non nulles. On peut alors conclure à l'aide de la fin du §3.

UNIVERSITÉ DE PARIS-SUD, ORSAY, FRANCE
*E-mail address*: clozel@math.u-psud.fr